\documentclass[11pt]{article}
\usepackage[rmargin=0.35in,right=0in, outer=0in,lmargin=0.35in,left=0in, inner=0in,tmargin=0.35in,bmargin=0.5in,paperwidth=5.5in,paperheight=8.5in]{geometry}
\usepackage{amsmath,amsxtra,amssymb,latexsym, amscd,amsthm}
\usepackage{mathptmx}
\usepackage{enumerate}  
\usepackage{setspace}%

\usepackage{mathtools, graphicx}
\usepackage{hyperref}
\usepackage[mathscr]{euscript}

\usepackage{nag, todonotes}

\newcommand{\R}{\mathbb{R}}

\newcommand{\Z}{\mathbb{Z}}

\newcommand{\F}{\mathcal{F}}
\newcommand{\T}{\mathcal{T}}

\newcommand{\G}{\mathcal{G}}
\newcommand{\h}{\mathcal{H}}
\newtheorem{theorem}{Theorem}[section]
\newtheorem{definition}{Definition}[section]
\newtheorem{corollary}{Corollary}[theorem]
\newtheorem{lemma}[theorem]{Lemma}
\theoremstyle{remark}
\newtheorem{remark}{Remark}[theorem]
\newcommand{\abs}[1]{\left\lvert #1 \right\rvert}

\begin{document}

\title{\textbf{{\Large A NEW CLASS OF GEOMETRICALLY DEFINED HYPERGRAPHS ARISING FROM THE HADWIGER-NELSON PROBLEM\thanks{This work was supported by NSF DMS grant no. 1950563}}}}

\author{
\small \textbf{Sean Fiscus}\\
\small Duke University\\
\small sean.fiscus@duke.edu \\ \\
\small \textbf{Eric Myzelev}\\
\small University of Pennsylvania\\ 
\small myzelev@sas.upenn.edu \\ \\
\small \textbf{Hongyi Zhang}\\
\small Haverford College \\
\small hzhang3@haverford.edu
}

\date{}
\maketitle 


\abstract{
There is a famous problem in geometric graph theory to find the chromatic number of the unit distance graph on Euclidean space; it remains unsolved. A theorem of Erdős and De-Bruijn simplifies this problem to finding the maximum chromatic number of a finite unit distance graph. Via a construction built on sequential finite graphs obtained from a generalization of this theorem, we have found a class of geometrically defined hypergraphs of arbitrarily large edge cardinality, whose proper colorings exactly coincide with the proper colorings of the unit distance graph on $\R^d$. We also provide partial generalizations of this result to arbitrary real normed vector spaces.
}

\section{Introduction}

$\Z^+$ will denote the set of positive integers. A hypergraph is a pair $\h=(V, E)$ in which $V$ is a non-empty set, the set of vertices of $\h$, and $E\subseteq 2^V$ (i.e. $E$ is a set of subsets of $V$) satisfying $e\in E$ implies $|e|\geq 2$. (In other definitions, singletons may be allowed in $E$ and/or $E$ may be a multi-set.) $E$ is the set of hyperedges, or just edges, of $\h$. A proper coloring of $\h$ is a function $\varphi:V\to C=$ some set of colors, such that no $e\in E$ is ``monochromatic with respect to $\varphi$". This means that for each $e\in E$, $\varphi|_e$ is not constant. The chromatic number of $\h$ is the smallest cardinality $|C|$ such that there is a proper coloring $\varphi:V\to C$.

If $\h=(V,E)$ is a hypergraph and $U\subseteq V$, $|U|\geq 2$, the subhypergraph of $\h$ induced by $U$ is $\h|_U=(U, E\cap 2^U)$. That is, the vertex set of $\h|_U$ is $U$ and the edges of $\h|_U$ are the edges of $\h$ that lie in $U$.
A hypergraph $\h=(V,E)$ in which $|e|=2$ for every $e\in E$ is a simple graph. The De Bruijn-Erdős theorem [1] concerning the chromatic numbers of infinite graphs is as follows:

If $m\in \Z^+, \h=(V, E)$ is a graph, and $\chi(\h|_F)\leq m$ for every finite set $F\subseteq V$, then $\chi(\h)\leq m$.

To prove our results we shall need a known generalization of this theorem to hypergraphs, which will be provided in the next section.

\section{De Bruijn-Erdős for hypergraphs with finite edges}

The following theorem generalizes the De Bruijn-Erdős Theorem to hypergraphs. A different proof may be found in \cite{soifer}, Chapter 26.
\begin{theorem}[D-E]
    Suppose that $m\in \Z^+ $, $\h=(V,E)$ is a hypergraph with $2\leq |e|<\infty$ for every $e\in E$, and $\chi(\h|_F)\leq m$ for every finite subset $F$ of $V$. Then $\chi(\h)\leq m$.
\end{theorem}
\begin{remark}
Since enlarging the edges of a hypergraph -- i.e. putting new vertices in some of the edges to join those that were there to begin with -- would seem to make it easier to avoid monochromatic edges when coloring the vertices, we are wondering if we really need the hypothesis, in Theorem D-E, that all edges of $\h$ are finite subsets of V. In the proof to come the hypothesis is actually used, and we do not see a way to avoid that use, except by trading the hypothesis for other, clumsier hypotheses. Is there a mathematical logician in the house?
\end{remark}
\begin{proof}
    Suppose that $\h=(V,E)$, and $m\in\Z^+$ satisfy the hypothesis of Theorem D-E. Let $[m]=\{1,\dots,m\}$. A coloring of V with colors $1,\dots,m$ is an element of $X=[m]^V=\prod_{v\in V}[m]$, the Cartesian product of $[m]$ with itself $|V|$ times. Let $[m]$ have the discrete topology, in which every singleton is an open set. Since $[m]$ is finite, obviously $[m]$ is compact (in this or any other topology). By the Tychonoff theorem \cite{munkres}, X, with the product topology, is compact.
    
    In the usual definition of the product topology on X, the basic open neighborhoods of a coloring $\varphi\in X$ are the sets $N(\psi, F)=\{\psi \in X |\  \psi|_F =\varphi|_F\}$,  $F\in \F(V)=\{$finite subsets of $ V\}$. 
    For $F\in\F(V)$, let $Y_F=\{\varphi\in X |\ \varphi|_F: F\rightarrow [m]$ is a proper coloring of $\h|_F\}$. 
    Since $\chi(\h|_F)\leq m, Y_F$ is nonempty.
    Also, $Y_F$ is closed in $X$: if $\psi\in X\backslash Y_F$, then $\psi|_F$ is not a proper coloring of $H_F$, so for some $e\in E\cap 2^F$, $\psi|_e$ is constant (i.e., the vertices in $e$ are all assigned the same $j\in[m]$ by $\psi$). 
    Any coloring in $N(\psi,e)$ will assign that same $j$ to each vertex of $e$; therefore, $N(\psi,e)\subseteq X\setminus Y_F$. Thus $X\backslash Y_F$ is open in the product topology, so $Y_F$ is closed.

    If $F_1,\dots, F_t\in \F(V)$, then $F_1\cup\dots\cup F_t\in \F(V)$, and we have \[\emptyset\neq Y_{\bigcup_{i=1}^tF_i} \subseteq Y_{F_1}\cap\dots\cap Y_{F_t}\] Thus $\{Y_F|F\text{ is in }\F(V)\}$ has the finite intersection property: any intersection of finitely many sets from the family is non-empty.

    In a compact topological space, the intersection of all the sets in a collection of closed sets with the finite intersection property is non-empty. Therefore $\bigcap_{F\in \F(V)} Y_F\neq \emptyset$. Suppose $\varphi \in \bigcap_{F\in \F(V)} Y_F$. We claim that $\varphi: V\rightarrow[m]$ is a proper coloring of $\h$. Suppose that $e\in E$. Then $e$ is a finite subset of $V$. Since $e$ is an edge of $\h\mid_e$, $\varphi\in Y_e$ implies that $\varphi$ assigns at least two colors from [m] to the elements of $e$.
\end{proof}

\section{Hypergraphs Equivalent to the Unit Distance Graph}

The motivation for this result is the following question: Does there exist a finite set of triangles $S$ in $\R^d$ such that the number of colors in a coloring of $\R^d$ required to forbid monochromatic copies of triangles in $S$ is the same as the chromatic number of the unit Euclidean distance graph on $\R^d$? The answer is yes. In fact, we prove a stronger result: in addition to generalizing triangles to arbitrary $m$-point sets, we also show that there is such a set $S$ so that a coloring $\varphi$ of $\R^d$ forbids congruent copies of $m$-gons in $S$ if and only if $\varphi$ forbids unit distance. To make this precise, we introduce a notion of equivalence of hypergraphs.

\begin{definition}
Let $S$ be a set, and $\h=(S, E_\h), \mathcal{G}=(S, E_\mathcal{G})$ be hypergraphs where $S$ is the vertex set. We say $\h$ is \underline{equivalent} to $\mathcal{G}$ if $\chi(\h)=\chi(\mathcal{G})$ and $\varphi:S\to C$, such that $|C|=\chi(\h)=\chi(\mathcal{G})$, is a proper coloring of $\h$ if and only if $\varphi$ is a proper coloring of $\mathcal{G}$.
\end{definition}

For the purposes of this paper, we define an $m$-gon as an arbitrary set of cardinality $m$. This means that some vertices of the $m$-gon can be collinear, which contradicts the standard geometric definition.

We also define a unit $m$-gon to be an $m$-gon in $\R^d$ as defined above with the following additional property: there exists at least one pair of points $x,y$ in $m$ which are Euclidean distance 1 apart. 

\begin{theorem}
Let $M\subset 2^{\R^d}$ be a non-empty finite set of $m$-gons, for some $m\geq 2$. Define $\h(M)$ as the hypergraph on $\R^d$ with edge set $E=\{X\subset\R^d\mid X \text{ is congruent in $\R^d$ to}$ $\text{some $T\in M$}\}$. 

Then there must exist some finite set $S$ of $(m+1)$-gons, such that $\h(M)$ is equivalent to $\h(S)$, where $\h(S)$ is an $(m+1)$-uniform hypergraph with edge set $E=\{Y\in\R^d\mid Y$ is congruent in $\R^d$ to some element of $S\}$.
\end{theorem}

\begin{proof}

We shall proceed by recursively obtaining sets $S_1,S_2,\dots$; $F_1, F_2, \dots$  satisfying:
    \begin{enumerate}
        \item $S_1\subseteq S_2\subseteq \dots$
        \item Each $S_j$ is a finite set of $(m+1)$-gons in $\R^d$, with each $(m+1)$-gon containing some $m$-gon $X\in M$ as a subset of its points.
        \item Defining $E_j:=\{$congruent copies in $\R^d$ of the $(m+1)$-gons in $S_j\}$ and $\h_j:=(\R^d,E_j)$, we obtain $F_j$, a finite subset of $\R^d\setminus M$ such that $\chi(\h_j|_{F_j})=\chi(\h_j)$.
        \item   $S_{j+1}=S_j\cup \{X\cup \{z\}|z\in F_j, X\in M\}$.
    \end{enumerate}

 Before giving the recursion, let us note that if the $S_j, F_j, \h_j=(\R^d, E_j)$ are as above then we have the following observations.
    \begin{enumerate}[(i)]
        \item Since, for each $j=1,2,\dots$ and $e\in E_j$, $e$ contains a copy of some $X\in M$, $\chi(\h_j)\leq \chi(\h(M))$.
        \item In view of the definition of $E_j$, $S_1\subseteq S_2\subseteq \dots$ implies that $E_1\subseteq E_2\subseteq \dots$, and thus $\chi(\h_1)\leq \chi(\h_2)\leq \dots$
        \item $\chi(\h(M))$ is finite. To see this, observe that in $\R^d$ with the Euclidean norm $||\  ||$, the sets $\{x,y\}$ congruent to a two-set $\{u,v\}$ are just the sets satisfying $||x-y||=||u-v||$. Therefore, a coloring of $\R^d$ forbidding congruent copies of $\{u,v\}$ is, in other jargon, the same as a coloring which forbids the distance $||u-v||$. For every positive distance, the smallest number of colors needed to forbid that distance is $\chi(\R^d,1)$, the chromatic number of the Euclidean unit distance graph on $\R^d$. 
\\Suppose that $|M|=n$ (recall that $M$ is finite). From each $T\in M$ select a $2$-set $D$; let the $2$-sets selected be \\$D_1,D_2,...,D_n$ and the distances determined by the $2$ points in these sets be $d_1,..,d_n$. For each $i=1,...,n$, let $\varphi_i$ be a coloring of $\R^n$ with $\chi(\R^d,1)$ colors that forbids the distance $d_i$. Now color $\R^d$ by an assignment $\psi$ of n-tuples: $\psi(r)=(\varphi_1(r),...,\varphi_n(r))$. We now have colored $\R^d$ with $\chi(R^d,1)^n<\infty$ colors, and it is easy to see that $\psi$ is a proper coloring of $\h(M)$: For any $T'\subseteq \R^d$ congruent to some $T\in M$, $\T'$ will contain a doubleton congruent to one of the $D_i$; to the two vectors in that doubleton, $\varphi_i$ will assign different colors, which means that $\psi$ will assign different n-tuples to the $2$ vectors, which means that $T'$ is not monochromatic.
        \item  By (i), (iii), and Theorem D-E, it follows that for each $j$ there is a finite set $F_j\in \R^d$ such that $\chi(\h_j|_{F_j})=\chi(\h_j)$.
        \item Suppose $\rho$ is an isometry of $\R^d$ and $F$ is a finite non-empty subset of $\R^d$. Because $\rho$ maps each $e\in E_j$ to an edge $\rho(e)\in E_j$, $\chi(\h_j|_F)=\chi(\h_j|_{\rho(F)})$.  \\
    \end{enumerate}

 The recursion: \\

    Let $S_1=\{X\cup\{a\}\mid X\in M\}$, for some $a\in\R^d\backslash\bigcup_{T\in M} T$. Let $F_1\subseteq \R^d$ be a finite set such that $\chi(\h_1|_{F_1})=\chi(\h_1)$; as explained in (iv) above, Theorem D-E guarantees the existence of such an $F_1$, and by (v), we can assume $F_1\cap (\bigcup_{T\in M} T)=\emptyset$; if $F_1\cap (\bigcup_{T\in M} T)\neq \emptyset$, replace $F_1$ by a translate of itself.

    From there, the recursion is dictated in $3$ and $4$: From $S_j$ we get $E_j$; the existence of $F_j$ is guaranteed. Then we define $S_{j+1}$ by 4, above, and roll on.
    
    Since the integer sequence $(\chi(\h_j))_j$ is non-decreasing and bounded above by $\chi(\h(M))$, clearly it will be eventually constant. If that eventual constant value were $\chi(\h(M))$, we would be almost done, except for showing that every proper coloring of $\h_k$ also properly colors $\h(M)$. All will be accomplished by the following.

    Clearly $\chi(\h_1)\geq 2>1$ whereas $\chi(\h_k)\leq \chi(\h(M))\leq k$ for $k\geq \chi(\h(M))$. Therefore there is a first value of $k\in \Z^+$ such that $\chi(\h_k) \leq k$.

    We have $k>1$ and \[k-1<\chi(\h_{k-1})\leq \chi(\h_k)\leq k\]  whence $\chi(\h_{k-1})=\chi(\h_k)=k\leq \chi(\h(M))$.

    Let $\varphi:\R^d\rightarrow \{1,\dots,k\}$ be a proper coloring of $\h_k$, and since $S_{k-1}\subset S_k$, $\varphi$ is also a proper coloring of $\h_{k-1}$. If $\varphi$ is a proper coloring of  $\h(M)$, then $\chi(\h(M))=k$, and, since $\varphi$ is an arbitrarily chosen proper coloring of $\h_k$, the claim of this theorem will be affirmed, with $S=S_k$.

    Suppose, on the contrary, that $\varphi$ does not properly color $\h(M)$, implying that for some $X'=\{a_1,...,a_m\}$ a congruent copy of some $X\in M$, $\varphi(a_1)= \varphi(a_2)=...=\varphi(a_m).$ We can, without loss of generality, convene that $\varphi(a_1)=...=\varphi(a_m)=k$. Let $\rho:\R^d\rightarrow\R^d$ be an isometry such that $\rho(X)=X'$. Consider any point $\rho(z)\in \rho(F_{k-1})$, and note that if it was colored with color $k$ then $X'\cup\{\rho(z)\}$ would be monochromatic and congruent to $X\cup\{z\}$, an edge in $\h_{k-1}$, which is impossible since $\varphi$ properly colors $\h_{k-1}$. Thus $\varphi$ is a proper coloring of $\h_{k-1}\mid_{F_{k-1}}$ with $k-1$ colors. This means that $\chi(\h_{k-1})=\chi(\h_{k-1}\mid_{F_{k-1}})\leq k-1$, and contradicts the fact that $\chi(\h_{k-1})=k$. $\qed$
\end{proof}

\begin{corollary}
For all integers $m\geq 2$, there exists a finite set $S$ of unit $m$-gons, such that $\h(S))$ is equivalent to $(\R^d,1)$, where $\h(S)=\h(\R^d,E)$ and $E=\{X\subset\R^d\mid X$ is congruent in $\R^d$ to some element of $S\}$ and $(\R^d,1)$ is the Euclidean unit distance graph on $\R^d$.
\end{corollary}

\begin{proof}

We prove the statement by induction. When $m=2$, the corollary is trivially true.

Now consider $m>2$ and suppose that the corollary holds true for $m-1$. Then there exists a finite set $S_{m-1}$ of $m-1$-gons as asserted in the Corollary. By the Theorem, there exists a finite set $S_m$ of $m$-gons such that $\chi(\h(S_m))=\chi(\h(S_{m-1}))$, and any proper coloring of $\h(S_m)$ with $\chi(\h(S_m)$ colors is also a proper coloring of $\h(S_{m-1})$. 

According to the inductive hypothesis, \\ $\chi(\h(S_{m-1}))=\chi(\R^d,1)$, and any proper coloring of $\h(S_{m-1})$ with $\chi(\h(S_{m-1})$ colors is also a proper coloring of the Euclidean unit distance graph on $\R^d$. Thus, $\chi(\h(S_m))=\chi(\R^d,1)$, and any proper coloring of $\h(S_m)$ with $\chi(\h(S_m))$ colors is also a proper coloring of the Euclidean unit distance graph on $\R^d$. 

Also note that because there are two points distance 1 apart in all of the $m-1$-gons in $S_{m-1}$, for every $m$-gon in $S_m$, there must also be two points distance 1 apart. That is, $S_m$ is a finite set of unit $m$-gons.
\end{proof}

\section{Generalizing to Non-Euclidean Norms on $\R^d$}

In the preceding sections, distance in $\R^d$ was provided by the Euclidean norm, hereinafter to be denoted as $||\cdot ||_2$. Some, but not all, of Theorem 3.1 and its corollary survives generalization to the setting of a finite-dimensional normed vector space over $\R$. Without loss of generality, the vector space will be $\R^d$ and the norm will be denoted as $||\cdot ||$. 

Two sets $X,Y\subseteq\R^d$ are \underline{\textit{congruent}} \underline{\textit{copies}} of each other in $(\R^d,||\cdot ||)$ if and only if one of them is the image of the other under a composition, in either order, of a surjective linear isometry of $(\R^d,||\cdot ||)$ and a translation. With this definition of congruence, we lose one of the support beams to our geometric intuition that may seem essential to the proof of Theorem 3.1: there can exist $u,v,x,y\in\R^d$ such that $||u-v||=||x-y||>0$ and yet $\{u,v\}$ and $\{x,y\}$ are not congruent. 

However, note: if $\{u,v\}$ and $\{x,y\}$ are congruent, then $||u-v||=||x-y||$. Consequently, if $\varphi$ is a coloring of $\R^d,||\cdot ||)$ which forbids a distance $a>0$, then $\abs{\varphi(e)}>1$ if $e\subseteq\R^d$ contains two points a distance $a$ apart. 

For $a>0$ let $\chi((\R^d,||\cdot ||),a)$ denote the smallest $\abs{C}$ such that some coloring $\varphi:\R^d\rightarrow C$ forbids the distance $a$. By the properties of norms, it is clear that $\chi((\R^d,||\cdot ||), a)=\chi((\R^d,||\cdot ||),1)$ for all $a>0$. Also very important for our generalizations:

\begin{lemma} For all $||\cdot ||$, $\chi((\R^d,||\cdot ||),1)<\infty$.
\end{lemma}
This is well known, but we will provide a proof outline in an Appendix. 

From Lemma 4.1 and the definition of congruence we obtain, as in section 3, and by the same argument, the following. 

\begin{lemma}
    Suppose that $\S$ is a finite collection of subsets of $\R^d$, each subset with at least 2 elements, and $\h(\S)=(\R^d,E(\S))$ is defined by $E(\S)=\{T\subseteq\R^d|T\text{ is congruent}$ $\textrm{ to some } S\in \S\}$. Then $\chi(\h(\S))<\infty$. 
\end{lemma}

\begin{theorem}
    For any norm $||\cdot ||$ on $\R^d$, Theorem 3.1 holds with the Euclidean norm replaced by $||\cdot ||$, provided the phrase ``congruent in $\R^d$" is replaced by ``congruent in $(\R^d,||\cdot ||)$."
\end{theorem}
In view of Lemmas 4.1 and 4.2, the proof is straightforward; follow the path of argument in section 3.

However, generalizing Corollary 3.1.1 looks to us like a lost cause. In the case $||\cdot ||=||\cdot ||_2$, for any two unit vectors $u,v$, there is a linear isometry of $(\R^d,||\cdot ||_2)$ that takes $u$ into $v$, whence $\{\underline{0},u\}$ and $\{\underline{0},v\}$ are congruent. Thus the graph $((\R^d,||\cdot ||),1)$ is the hypergraph $\h(\{\underline{0},u\})$, which makes the base of the induction proof of the Corollary, when $m=2$, ``trivial".

Given a non-Euclidean norm $||\cdot ||$ on $\R^d$, we can, with reference to Theorem 4.3, find $\S_2\subseteq \S_3 \subseteq \dots$ such that $\S_m$ is a finite set of unit $m$-gons in $(\R^d,||\cdot ||)$, $m=2,3,\dots$, and the hypergraphs $\h(\S_m)$ are all equivalent, but we can only be sure that $\chi(\h(\S_m))\leq \chi((\R^d,||\cdot ||),1)$, and even if we happen to have equality, we see no way of assuring that every proper coloring of $\chi(\S_m)$ with number of colors equal to the chromatic number of graph will also properly color the unit distance graph on $(\R^d,||\cdot ||)$.

On the other hand, we have no proof that there is no generalization of Corollary 3.1.1 for some or all non-Euclidean norms $||\cdot ||$ on $\R^d$, $d>1$. Perhaps the question is worth investigating. (Yes, we are aware that Corollary 3.1.1 will hold for every norm arising from an inner product on $\R^d$; $\R^d$ with such a norm is linearly and isometrically isomorphic to $(\R^d,||\cdot ||_2)$.)

We can obtain $m$-uniform hypergraphs in $(\R^d,||\cdot ||)$ with the same chromatic number as $((\R^d,||\cdot ||),1)$ if we abandon the practice of defining edge sets as the sets congruent to one of a finite set of $m$-gons. However, our method does not provably produce hypergraphs equivalent to the unit distance graph on $(\R^d,||\cdot ||)$.

Let $Q$ be a collection of finite subsets of $\R^d$, each with at least 2 elements, and let 
$$E_0=\{e\subseteq \R^d\mid \textrm{for some $f\in Q$, $e$ and $f$ are congruent}\}$$
For each $t\in\Z^+$, let 
$$B_t=\{f\cup T\mid  f\in Q, T\subseteq \R^d, f\cap T=\emptyset, \textrm{ and } |T|=t \}$$
$$E_t=\{e\subseteq \R^d\mid\text{ for some } g\in B_t\text{, }e\text{ and } g \text{ are congruent}\}$$
It will be important to notice that because, for each $f\in Q$, we form infinitely many sets in $B_t$ by taking the union of $f$ with each and every $t-$subset of $\R^d\setminus f$, it follows that if $e\in E_0$, $T'\subseteq \R^d\setminus e$, and $\abs{T'}=t$, then $e\cup T'\in E_t$. 

\begin{theorem}
    With $t\in\Z^+,Q,B_t,E_0,$ and $E_t$ as above, let $\h_0=(\R^d,E_0)$ and $\h_t=(\R^d,E_t)$. If $\chi(\h_0)<\infty$ then $\chi(\h_t)=\chi(\h_0)$. 
\end{theorem}

\begin{proof}
        Since each $e\in E_t$ contains an $e'\in E_0$, it follows that a proper coloring of $\h_0$ will also serve as a proper coloring of $\h_t$. Thus, \[\chi(\h_t)\leq\chi(\h_0)<\infty\] By the same argument, \[\chi(\h_t)\leq\chi(\h_{t-1})\leq\chi(\h_0)\] It suffices to show that $\chi(\h_t)=\chi(\h_{t-1})$ for each $t\in\Z^+$. 

        Let $k=\chi(\h_t)$ and suppose that $k<\chi(\h_{t-1})$. Let $\varphi:\R^d\rightarrow\{1,...,k\}$ be such that no edge $g\in E_t$ is monochromatic, with reference to $\varphi$.

        Since $k<\chi(\h_{t-1})$, there must exist some $B\in E_{t-1}$ such that $\varphi$ assigns the same color to every element of $B$. Without loss of generality, we can assume that this color is $k$. 

        For every $u\in\R^d\setminus B$, $B\cup\{u\}\in E_t$. Therefore \\ $\varphi(u)\in\{1,...,k-1\}$; otherwise, $B\cup\{u\}$ would be monochromatic under the coloring $\varphi$. Thus, $\varphi$ restricted to $\R^d\setminus B$ is a proper coloring of $\h_t\mid_{\R^d\setminus B}$ with $k-1$ colors. 

        We shall now use Theorem $D-E$ to prove the existence of a proper coloring of $\h_t$ with colors $\{1,...,k-1\}$ which will contradict $k=\chi(\h_t)$. Since this contradiction descends from the assumption that \\ $k=\chi(\h_t)<\chi(\h_{t-1})$ it will follow that $\chi(\h_t)=\chi(\h_{t-1})$ and the theorem will be proven. 

        Suppose $F\subset\R^d$ is finite and $\chi(\h_t\mid_F)=\chi(\h_t)$. We aim to show that $\h_t\mid_F$ can be properly colored with no more than $k-1$ colors, which will imply that \[\chi(\h_t)=\chi(\h_t|_F)\leq k-1<k=\chi(\h_t).\]
        
        Let $v\in \R^d$ be such that $(v+F)\cap B=\emptyset$. Then $\varphi$ colors $v+F$ with no more than $k-1$ colors so that for each $\alpha\in2^{v+F}\cap E_t$, $\varphi$ assigns more than one color to the elements of $\alpha$. Now color $F$ as follows: color $f\in F$ with $\varphi(v+f)$. Since every translate of every $\alpha\in E_t$ is in $E_t$, and no $\alpha\in 2^{v+F}\cap E_t$ is monochromatic under coloring by $\varphi$, we have what we wanted, a proper coloring of $\h_t\mid_F$ with colors from $\{1,...,k-1\}$. 
\end{proof}

\begin{corollary}
Let $||\cdot ||$ be a norm on $\R^d$. For each $m\in \Z^+$ such that $m>2$, define $E_m=\{T\subseteq \R^d \mid |T|=m$ and $T$ contains 2 points $ ||\cdot ||$-distance 1 apart$\}$. Let $\G_m=(\R^d, E_m)$. Then $\chi(\G_m)=\chi((\R^d,||\cdot ||), 1)$ for all $m$.
\end{corollary}

\begin{proof}
If $Q=\{\{\underline{0}, u\}|u\in\R^d\text{ and }||u||=1\}$ , then, in terms used in the Theorem's statement, $\h_0=((\R^d,||\cdot ||),1)$, the unit distance graph on $(\R^d,||\cdot ||)$ and, for each $m>2$, $\h_{t-2}=\G_m$. The conclusion follows from the Theorem. 
\end{proof}

\section{Explicit Construction of Hypergraph}

In this section, we are in $\R^d$ with the usual Euclidean norm. We will give a ``constructive" proof of a weaker version of Theorem 3.1.

By applying Theorem D-E for hypergraphs to Theorem 3.1, it can be shown that for every finite $m$-uniform hypergraph in $\R^d$, there must exist a finite $(m+1)$-uniform hypergraph of equal chromatic number. However, this deduction is \\non-constructive because all known proofs of Theorem D-E use the axiom of choice, and thus we have no control over what the finite hypergraph might look like. In this section, we show how to take a finite $m$-uniform hypergraph and construct a finite $(m+1)$-uniform hypergraph of equal chromatic number.   This allows us to ``construct", for any $m\in \Z^+$, $m>2$, a finite $m$-uniform hypergraph in $\R^d$ with chromatic number equal to $\chi(\R^d, 1)$. We use ``construct" in quotations, however, because this hypergraph construction will use a finite unit distance graph $G$ with chromatic number $\chi(\R^d, 1)$.

\begin{theorem}

Let $\h$ be a finite, $m$-uniform hypergraph with vertices in $\R^d$. There exists a finite, $(m+1)$-uniform hypergraph $\h'$ with vertices in $\R^d$ such that $\chi(\h')=\chi(\h)$.

\end{theorem}

\begin{proof}

Let $k=\chi(\h)$, and $F$ be the vertex set $V(\h)$. Note that for all translates $F+v\subset\R^d$ of $F$ where $v\in\R^d$, the chromatic number of corresponding hypergraph $\h+v$, with vertex set $V(\h+v)=F+v$ and $E(\h+v)=\{e+v\mid e\in E(\h)\}$, must also be $k$. 

We now construct an $(m+1)$-uniform hypergraph $\h'$ in the following way:

\begin{enumerate}
    \item Let the vertex set of $\h'$ be the union of $k$ disjoint translates of $F=V(\h)$, which we will call $F_1,F_2,...,F_k$. Denote by $\h_i$ the translate of $\h$ onto $F_i$. 
    \item Define the edge set of $\h'$ by \[E(\h')=\{\{e\cup v\}\mid e\in\h_i\text{, } v\in F_j\text{, where }j>i\}\]
In other words, $E(\h')$ consists of all $(m+1)$-gons such that $m$ points are from an edge in $\h_i$ and the remaining point is from some $F_j$ with $j>i$.
\end{enumerate}

By construction, $\h'$ is finite and $(m+1)$-uniform. Also, any proper coloring $\varphi$ of $\h$ can be extended to a proper coloring of $\h'$ by coloring each copy $F_i$ of $V(\h)$ in $V(\h')$ by the same colors assigned to $V(\h)$. Doing so, all $e\in E(\h_i)$ will be properly colored, implying that all  $e\in E(\h')$ must be properly colored by construction. Thus, we have shown that $k\geq\chi(H')$, and need only to confirm that $k\not>\chi(H')$. To show this, we demonstrate that any coloring of $\h'$ with fewer than $k$ colors must yield a monochromatic edge. 

Consider a coloring $\varphi:\R^d\to\{1,...,k-1\}$. It is clear $\varphi$ must not properly color $\h_i$ for each $\h_i\in \h_1,...,\h_k$. Thus for each $\h_i$, there must be some edge $e=\{a_1,...,a_m\}\in E(\h_i)$ such that $\varphi(a_1)=...=\varphi(a_m)$. 

We now show that there must be an edge $e\in\h'$ which $\varphi$ monochromatically colors. 

Let $e_1$ be a monochromatic edge in $E(\h_1)$. Without loss of generality suppose that $\varphi(e_1)=\{k-1\}$. If any vertex $v\in F_2\cup\ ...\ \cup F_k$ were colored by $k-1$, we would have found a monochromatic edge in $\h'$, by our construction of the edges in $\h'$. Thus, we can assume that $\varphi(v)\in\{1,...,k-2\}$ for $v\in F_2\cup\ ...\ \cup F_k$. 

Likewise, there must be a monochromatic edge $e_2$ in $E(\h_2)$, because $F_2$ is colored now with even fewer colors than $F_1$ was. Because $\varphi(F_2)\in\{1,...,k-2\}$, we can without loss of generality suppose that $\varphi(e_2)=\{k-2\}$. 

Continuing in this fashion until reaching $F_k$, we see that for each color $\alpha\in\{1,...,k-1\}$, there must appear a monochromatic edge in some $\h_i$ for $i\in\{1,...,k-1\}$. Now consider $\varphi(v)$ for $v\in F_k$. No matter which color is chosen for $v$, there will appear a monochromatic edge in $\h'$ of the form $e_i\cup v$ for some $e_i\in\h_i$. We see that it is impossible to properly color $\h'$ with fewer than $k$ colors, implying that $\chi(\h')=k$.  

\end{proof}

\begin{corollary}
    For all integers $m\geq 2$, there exists a finite, $m$-uniform hypergraph $\h'$ with vertices in $\R^d$ such that $\chi(\h)=\chi(\R^d,1)$.   
\end{corollary}

\begin{proof}

We prove the statement by induction. When $m=2$, the result follows directly from the D-E Theorem.

Now consider $m>2$ and suppose the theorem holds true for $m-1$. Then there exists a finite, $(m-1)$-uniform hypergraph $\h$ with vertices in $\R^d$ such that $\chi(\h)=\chi(\R^d,1)$. By the previous theorem, there must exist an $m$-uniform hypergraph $\h'$ with vertices in $\R^d$ such that $\chi(\h')=\chi(\h)=\chi(\R^d,1)$. $\qed$

\end{proof}
\begin{remark}
    When $m=2$, the hypergraph is not obtained constructively, but with this method, we still have control over the large-scale geometric structure of the $m$-uniform hypergraph based on how we position the finite hypergraphs when $m>2$, and hence the term  ``construct". 
\end{remark}

\section*{Appendix}

Summary of the proof that for any norm $||\cdot ||$ on $\R^d$, $\chi(\R^d,1)<\infty$. 

Let $||\cdot ||$ be a norm on $\R^d$ and let $||\cdot ||_\infty$ be the norm on $\R^d$ defined by $||(a_1,...,a_d)||_\infty$ $=\max{[\abs{a_i};1\leq i\leq d]}$.

The proof of this well-known result will rely on the even better-known fact that any two norms on $\R^d$ are \underline{equivalent}. The equivalence of $||\cdot||$ and $||\cdot ||_\infty$ means that there exist $c,C>0$ such that for all $u\in\R^d$, \[c||u||_\infty\leq||u||\leq C||u||_\infty.\] The rightmost inequality can be used to prove the existence of $\varepsilon>0$ such that the $||\cdot ||$-diameter of the d-dimensional cube $[0,\varepsilon]^d$ is $<1$. [Just take $0<\varepsilon<\frac{1}{c}.]$ Then the leftmost inequality above implies the existence of an integer $m>0$ such that for each of the unit coordinate vectors $e_j=(S_{1j},...,S_{dj})$ (in which $S_{ij}$ is the Kronecker delta), $||(m-1)e_j||>1$. 

Next we partition the cube $Q=[0,m\varepsilon)^d$ into $m^d$ little cubes, $\sum^{d-1}_{t=0}[j_t\varepsilon,(j_t+1)\varepsilon),$ $(j_0,...,j_{d-1})\in\{0,1,...,m-1\}^d$, and color $Q$ with $m^d$ colors; one to each little cube. Finally we color $\R^d$ with those same colors by tiling $\R^d$ with translates (with colors attached) of $Q$. 

This will be a proper coloring of $((\R^d,||\cdot ||),1)$ because any two points bearing the same color are either within the same little cube— no two points of which are  $||\cdot ||$-distance $\geq 1$ apart,— or in two different little cubes, in which case the distance between them is $>1$. 

Thus $\chi(\R^d,||\cdot ||)\leq m^d<\infty$. 

\section*{Acknowledgement}

We would like to acknowledge that this paper was written with much advice from Dr. Peter Johnson, our mentor during the 2023 NSF-funded summer REU program at Auburn University. We extend our deepest gratitude for his generous guidance on this project and others throughout our time at Auburn.

\end{document}